\documentclass[12pt]{report}
\usepackage{amssymb}
\usepackage{amsmath}

\begin{document}

\vspace{.2in}\parindent=0mm
\begin{center}

 {\bf\large {Gabor Frames on Local Fields of Positive Characteristic  }}

\parindent=0mm \vspace{.5in}
{\bf{ Firdous A. Shah$^{*}$ }}

\end{center}

\parindent=0mm \vspace{.1in}
{{\it\small$^{*}${Department of  Mathematics,  University of Kashmir, South Campus, Anantnag-192 101, Jammu and Kashmir, India. E-mail: $\text{fashah79@gmail.com}$}}

{\small{
\parindent=0mm \vspace{.2in}
{\bf{Abstract:}} Gabor frames have gained considerable popularity during the past decade, primarily due to their substantiated applications in diverse and widespread fields of engineering and science. Finding general and verifiable conditions which imply that the Gabor systems are Gabor frames is among the core problems in time-frequency analysis. In this paper, we give some simple and sufficient conditions that ensure a Gabor system $\left\{M_{u(m)b}T_{u(n)a}g=:\chi_{m}(bx)g\big(x-u(n)a\big)\right\}_{ m,n\in\mathbb N_0}$ to be a frame for  $L^2(K)$. The conditions proposed are stated in terms of the Fourier transforms of the Gabor system's generating functions.  

\parindent=0mm \vspace{.2in}

\parindent=0mm \vspace{.1in}
 }}

\parindent=0mm \vspace{.2in}
{\bf{1. Introduction}}

\parindent=0mm \vspace{.2in}

The notion of frame was first introduced by Duffin and Schaeffer  [1] in connection with some deep problems in non-harmonic Fourier series. Frames are basis-like systems that span a vector space but allow for linear dependency, which can be used to reduce noise, find sparse representations, or obtain other desirable features unavailable with orthonormal bases. The idea of Duffin and Schaeffer did not generate much interest outside non-harmonic Fourier series until the seminal work by Daubechies, Grossmann, and Meyer [2]. They combined the theory of continuous wavelet transforms with the theory of frames to introduce wavelet (affine) frames for $L^2(\mathbb R)$.  After their work, the theory of frames began to be studied widely and deeply. Today, the theory of frames has become an interesting and fruitful field of mathematics with abundant applications in signal processing, image processing, harmonic analysis, Banach space theory, sampling theory,  wireless sensor networks, optics, filter banks, quantum computing, and medicine and so on. An introduction to the frame theory and its applications can be found in [3, 4].

\parindent=8mm \vspace{.2in}
Gabor frames form a special kind of frames for $L^2(\mathbb R)$ whose elements are generated by time-frequency shifts of a single window-function or atom. More specifically, let $g\in L^2(\mathbb R)$ and $a,b\in \mathbb R^{+}$, we use ${\cal G}(g, a, b)$ to denote the {\it Gabor family} or {\it system} $\left\{ M_{mb}T_{na}g:m,n\in \mathbb Z\right\}$ generated by $g$ where $T_{na} f(x)=f(x-na)$ is the translation unitary operator and $M_{mb}f(x)=e^{2\pi i mb x}f(x)$ is the modulation unitary operator. The composition $M_{mb}T_{na}$ is called the time-frequency shift operator. The system $\left\{ M_{mb}T_{na}g:m,n\in \mathbb Z\right\}$ is called a {\it Gabor frame} if there exist constants $A,B>0$ such that

$$A\big\|f\big\|_{2}^{2}\le \sum_{m\in\mathbb Z}\sum_{n\in\mathbb Z}\left|\big \langle f, M_{mb}T_{na}g\big\rangle \right|^2\le B \big\|f\big\|_{2}^{2}, \quad \text{for all}~f\in L^2(\mathbb R).\eqno(1.1)$$

\parindent=8mm \vspace{.1in}
Gabor systems that form frames for $L^2(\mathbb R)$ have a wide variety of applications. In practice, once the window function has been chosen, the first question to investigate for Gabor analysis is to find the values of the time-frequency parameters $a, b$ such that $\left\{ M_{mb}T_{na}g\right\}_{m,n\in \mathbb Z}$ is a frame. A useful tool in this context is the Ron and Shen [5] criterion. By using this criterion, Gr\"{o}chenig et al.[6] have proved that the system $\left\{ M_{mb}T_{na}g\right\}_{m,n\in \mathbb Z}$ cannot be a frame for $a > 0$ and $b$ integer greater than $1$. Many results in this area, including necessary conditions and sufficient conditions have been established during the last two decades [7--10]. We refer the reader to the books [11, 12] for a comprehensive treatment of Gabor frames.

\parindent=8mm\vspace{.2in}
A field $K$ equipped with a topology is called a local field if both the additive and multiplicative groups of $ K$ are locally compact Abelian groups.  For example, any field endowed with the discrete topology is a local field. For this reason we consider only non-discrete fields. The local fields are essentially of two types (excluding the connected local fields $\mathbb R$ and $\mathbb C$). The local fields of characteristic zero include the $p$-adic field $\mathbb Q_p$. Examples of local fields of positive characteristic are the Cantor dyadic group and the Vilenkin $p$-groups. Local fields have attracted the attention
of several mathematicians, and have found innumerable applications not only in the number theory, but also in the representation theory, division algebras, quadratic forms and algebraic geometry. As a result, local fields are now consolidated as a part of the standard repertoire of contemporary mathematics. For
more details we referr to [13].

\parindent=8mm \vspace{.2in}

The local field $K$ is a natural model for the structure of Gabor frame systems, as well as a domain upon which one can construct Gabor basis functions. There is a substantial body of work that has been concerned with the construction of Gabor frames on $K$, or more generally, on local fields of positive characteristic. Jiang et al.[14] constructed Gabor frames on  local fields of positive characteristic using  basic concepts of operator theory and have established a necessary and  sufficient conditions for the system  $\left\{ M_{u(m)b}T_{u(n)a}g=:\chi_{m}(bx)g\big(x-u(n)a\big)\right\}_{ m,n\in\mathbb N_0}$ to be a frame for  $L^2(K)$. Recently, Shah [15] established a complete characterization of Gabor frames on local fields by virtue of two basic equations in the Fourier domain and show how to construct an orthonormal Gabor basis for $L^2(K)$. Recent results related to wavelet and Gabor frames on local fields of prime characteristic can be found in [16--20] and the references therein.

\parindent=8mm \vspace{.2in}
In this article, we continue our investigation on Gabor frames on local fields and will present generalized inequalities for Gabor frames on local fields of positive characteristic via Fourier transform. The inequalities we proposed are stated in terms of the Fourier transforms of the Gabor system's generating functions, and the inequalities are better than that of Li and Jiang [14].  Although we consider a one-dimensional case here, our results are easily generalized to multi-dimensional Gabor systems on local fields of positive characteristic.

\pagestyle{myheadings}

\parindent=8mm \vspace{.2in}
The paper is organized as follows. In Section 2, we discuss some preliminary facts about local fields of positive characteristic and state the main results. Section 3 gives the proofs of the results.

\parindent=0mm \vspace{.2in}

{\bf{2.  Preliminaries on Local Fields }}

\parindent=0mm \vspace{.2in}
Let $K$ be a field and a topological space. Then $K$ is called a {\it local field} if both  $K^+$ and $K^*$ are locally compact Abelian groups, where  $K^+$ and $K^*$ denote the additive and multiplicative groups of $K$, respectively. If $K$ is any field and is endowed with the discrete topology, then $K$ is a local field. Further, if $K$ is connected, then $K$ is either $\mathbb R$ or $\mathbb C$. If $K$ is not connected, then it is totally disconnected. Hence by a local field, we mean a field $K$ which is locally compact, non-discrete and totally disconnected. We use the notation of the book by Taibleson [13]. Proofs of all the results stated in this section can be found in the book [13].

\parindent=8mm \vspace{.2in}
Let $K$ be a local field. Let $dx$ be the Haar measure on the locally compact Abelian group $K^{+}$. If $\alpha\in K$ and $\alpha\ne 0$, then $d(\alpha x)$ is also a Haar measure. Let $d(\alpha x)=|\alpha|dx$. We call $|\alpha|$ the {\it absolute value} of $\alpha$. Moreover, the  map $x\to |x|$ has the following properties: (a) $|x|=0$ if and only if $x = 0;$ (b) $|xy|=|x||y|$ for all $x, y \in K$; and (c) $|x+y|\le \max \left\{ |x|, |y|\right\}$ for all $x, y\in K$. Property (c) is called the {\it ultrametric inequality.} The set ${\mathfrak D}= \left\{x \in K: |x| \le 1\right\}$ is called the {\it ring of integers} in $K.$ Define ${\mathfrak B}= \left\{x \in K: |x| < 1\right\}$. The set ${\mathfrak B}$ is called the {\it prime ideal} in $K$. The prime ideal in $K$ is the unique maximal ideal in ${\mathfrak D}$ and hence as result ${\mathfrak B}$ is both principal and prime. Since the local field $K$ is totally disconnected, so there exist an element of ${\mathfrak B}$ of maximal absolute value. Let $\mathfrak p$ be a fixed element of maximum absolute value in ${\mathfrak B}$. Such an element is called a {\it prime element} of $K.$ Therefore, for such an ideal ${\mathfrak B}$ in ${\mathfrak D}$, we have ${\mathfrak B}= \langle \mathfrak p \rangle=\mathfrak p {\mathfrak D}.$ As it was proved in [13], the set ${\mathfrak D}$ is compact and open. Hence, ${\mathfrak B}$ is compact and open. Therefore, the residue space ${\mathfrak D}/{\mathfrak B}$ is isomorphic to a finite field $GF(q)$, where $q = p^{k}$ for some prime $p$ and $k\in\mathbb N$.

\parindent=8mm \vspace{.2in}
Let ${\mathfrak D}^*= {\mathfrak D}\setminus {\mathfrak B }=\left\{x\in K: |x|=1   \right\}$. Then, it can be proved that ${\mathfrak D}^*$ is a group of units in $K^*$ and if $x\not=0$, then we may write $x=\mathfrak p^k x^\prime, x^\prime\in {\mathfrak D}^*.$ For a proof of this fact we refer to [13]. Moreover, each ${\mathfrak B}^k= \mathfrak p^k {\mathfrak D}=\left\{x \in K: |x| < q^{-k}\right\}$ is a compact subgroup of $K^+$ and usually  known as the {\it fractional ideals} of $K^+$. Let ${\cal U}= \left\{a_i \right\}_{ i=0}^{q-1}$ be any fixed full set of coset representatives of ${\mathfrak B}$ in ${\mathfrak D}$, then every element $x\in K$ can be expressed uniquely  as $x=\sum_{\ell=k}^{\infty} c_\ell \mathfrak p^\ell $ with $c_\ell \in {\cal U}.$ Let $\chi$ be a fixed character on $K^+$ that is trivial on ${\mathfrak D}$ but is non-trivial on  ${\mathfrak B}^{-1}$. Therefore, $\chi$ is constant on cosets of ${\mathfrak D}$ so if $y \in {\mathfrak B}^k$, then $\chi_y(x)=\chi(yx), x\in K.$ Suppose that $\chi_u$ is any character on $K^+$, then clearly the restriction $\chi_u|{\mathfrak D}$ is also a character on ${\mathfrak D}$. Therefore, if $\left\{u(n): n\in\mathbb N_0\right\}$ is a complete list of distinct coset representative of ${\mathfrak D}$ in $K^+$, then, as it was proved in [13], the set  $\left\{\chi_{u(n)}: n\in\mathbb N_0\right\}$   of distinct characters on ${\mathfrak D}$ is a complete orthonormal system on ${\mathfrak D}$.

\parindent=8mm \vspace{.2in}
The Fourier transform $\hat f$ of a function $f \in L^1(K)\cap L^2(K)$ is defined by\\
$$\hat f(\xi)= \displaystyle \int_K f(x)\overline{ \chi_\xi(x)}dx.\eqno(2.1)$$
It is noted that\\
$$\hat f(\xi)= \displaystyle \int_K f(x)\,\overline{ \chi_\xi(x)}dx= \displaystyle \int_K f(x)\chi(-\xi x)dx.$$

\parindent=0mm \vspace{.2in}
Furthermore, the properties of Fourier transform on local field $K$ are much similar to those of on the real line. In particular Fourier transform is unitary on $L^2(K)$.

\parindent=8mm \vspace{.2in}
We now impose a natural order on the sequence $\{u(n)\}_{n=0}^\infty$. We have ${\mathfrak D}/ \mathfrak B \cong GF(q) $ where $GF(q)$ is a $c$-dimensional vector space over the field $GF(p)$. We choose a set $\left\{1=\zeta_0,\zeta_1,\zeta_2,\dots,\zeta_{c-1}\right\}\subset {\mathfrak D^*}$ such that span  $\left\{\zeta_j\right\}_{j=0}^{c-1}\cong GF(q)$. For $n \in \mathbb N_0$ satisfying
$$0\leq n<q,~~n=a_0+a_1p+\dots+a_{c-1}p^{c-1},~~0\leq a_k<p,~~\text{and}~k=0,1,\dots,c-1,$$

\parindent=0mm \vspace{.1in}
we define
$$u(n)=\left(a_0+a_1\zeta_1+\dots+a_{c-1}\zeta_{c-1}\right){\mathfrak p}^{-1}.\eqno(2.2)$$

\parindent=0mm \vspace{.1in}
Also, for $n=b_0+b_1q+b_2q^2+\dots+b_sq^s, ~n\in \mathbb N_{0},~0\leq b_k<q,k=0,1,2,\dots,s$, we set

$$u(n)=u(b_0)+u(b_1){\mathfrak p}^{-1}+\dots+u(b_s){\mathfrak p}^{-s}.\eqno(2.3)$$

\parindent=0mm \vspace{.1in}
This defines $u(n)$ for all $n\in \mathbb N_{0}$. In general, it is not true that $u(m + n)=u(m)+u(n)$. But, if $r,k\in\mathbb N_{0}\; \text{and}\;0\le s<q^k$, then $u(rq^k+s)=u(r){\mathfrak p}^{-k}+u(s).$ Further, it is also easy to verify that $u(n)=0$ if and only if $n=0$ and $\{u(\ell)+u(k):k \in \mathbb N_0\}=\{u(k):k \in \mathbb N_0\}$ for a fixed $\ell \in \mathbb N_0.$ Hereafter we use the notation $\chi_n=\chi_{u(n)}, \, n\ge 0$.

\parindent=8mm \vspace{.2in}
Let the local field $K$ be of characteristic $p>0$ and $\zeta_0,\zeta_1,\zeta_2,\dots,\zeta_{c-1}$ be as above. We define a character $\chi$ on $K$ as follows:
$$\chi(\zeta_\mu {\mathfrak p}^{-j})= \left\{
\begin{array}{lcl}
\exp(2\pi i/p),&&\mu=0\;\text{and}\;j=1,\\
1,&&\mu=1,\dots,c-1\;\text{or}\;j \neq 1.
\end{array}
\right. \eqno(2.4)
$$

\parindent=8mm \vspace{.1in}
We also denote the test function space on $K$ by $\Omega$, i.e., each function $f$ in $\Omega$ is a finite linear combination of functions of the form ${\bf 1}_k(x-h), h\in  K, k\in\mathbb Z$, where ${\bf 1}_k$ is the characteristic function of ${\mathfrak B}^k$. Then, it is clear that $\Omega$ is dense in $L^p( K),  1\le p <\infty$, and each function in $\Omega$ is of compact support and so is its Fourier transform.

\parindent=8mm \vspace{.2in}
For a given $\psi \in L^2(K)$, define the  Gabor system

$${\cal G}(g, a, b):= \Big\{M_{u(m)b}T_{u(n)a}g=:\chi_{m}(bx) g\big(x-u(n)a\big): n,m\in\mathbb N_0\Big\}.\eqno(2.5)$$

\parindent=0mm \vspace{.1in}
We call the Gabor  system ${\cal G}(g, a, b)$ a Gabor frame for $L^2(K)$, if there exist positive numbers $0 < C \le D < \infty$ such that for all $f\in  L^2(K)$

$$C\big\|f\big\|_{2}^{2}\le \sum_{m\in\mathbb N_{0}}\sum_{n\in \mathbb N_0} \left|\big\langle f, M_{u(m)b}T_{u(n)a}g\big\rangle\right|^2 \le D\big\|f\big\|_{2}^{2}.\eqno(2.6)$$

\parindent=0mm \vspace{.0in}
 Before stating our results, we introduce some notations. For any $g\in L^2(K)$ and $a, b>0$. We set
$$\begin{array}{rcl}
\Delta_k(\xi)&=&\displaystyle\sum_{m \in \mathbb{N}_0}\Big| \hat{g}\big(\xi -bu(m) \big) \hat{g}\big( \xi-bu(m) +a^{-1}u(k)\big )\Big|,\qquad\\\\
\alpha_k &=&\displaystyle \text{ess}\, \sup_{\xi}  \Delta_k(\xi),~k \in \mathbb N_{0},\quad \beta=\sum_{k \in \mathbb{N}} \alpha_{k},\quad \gamma=\text{ess} \inf_{\xi} \Delta_{0}(\xi),\qquad\\\\
\Lambda_k(\xi)&=&\displaystyle\sum_{k \in \mathbb{N}_0}\hat{g}\big(\xi-bu(m) \big) \hat{g}\big(\xi-bu(m)+a^{-1}u(k) \big),\qquad\\\\
\delta_k &=&\displaystyle\text{ess}\sup_{\xi} \big |\Lambda_k(\xi) \big |, ~k \in \mathbb N_{0},\quad \mu=\sum_{k \in \mathbb{N}}\delta_k\,.
\end{array}$$

\parindent=0mm \vspace{.0in}

Based on these notations, the authors in [14] established the following result.

\parindent=0mm \vspace{.2in}
{\bf{Theorem 2.1.}} {\it Let $a, b>0$ and $g\in L^2(K)$. If $\alpha_{0}, \beta$ and $\gamma$ satisfy}

$$\beta< \gamma\le \alpha_0 < \infty, $$

\parindent=0mm \vspace{.1in}
{\it then system  in (2.5)  constitutes a Gabor frame for $L^2(K)$ with bounds $\dfrac{C_{1}}{a}$ and $\dfrac{D_{1}}{a}$, where $C_{1}=\gamma-\beta$ and $D_{1}= \alpha_{0}+\beta.$ }

\parindent=0mm \vspace{.2in}
Motivating by the fundament works in [14,15], we will give two new sufficient conditions of Gabor frame on local fields of prime characteristic in this paper. The conditions obtained are better than that of one in Theorem 2.1.

\parindent=8mm \vspace{.1in}
Now, the first result of the paper is stated as follows.

\parindent=0mm \vspace{.2in}
{\bf{Theorem 2.2.}} {\it Let $a, b>0$ and $g\in L^2(K)$. If $\alpha_{0}, \gamma$ and $\mu$ satisfy}

$$\mu < \gamma\le \alpha_0 < \infty,\eqno(2.7) $$

\parindent=0mm \vspace{.1in}
{\it then the Gabor system ${\cal G}(g, a, b)$ as defined in (2.5) is a frame for $L^2(K)$ with bounds $\dfrac{C_{2}}{a}$ and $\dfrac{D_{2}}{a}$, where $C_{2}=\gamma-\mu$ and $D_{2}= \alpha_{0}+\mu.$ }

\parindent=0mm \vspace{.2in}
{\it Remark 1.} It is easy to see that $\mu\le \beta$, so the frame bounds in Theorem 2.2 are better than ones in Theorem 2.1.

\parindent=8mm \vspace{.2in}
Next, we prove a more general result which includes not only the results of Theorem 2.1 and  2.2 as special cases, but also leads to a standard development of
interesting generalizations of some well-known results.  To do so, we set

$$\sigma=\text{ess}\sup_{\xi}\sum_{k\in \mathbb N}\big|\Lambda_{k}\big|. $$

\parindent=0mm \vspace{.1in}
{\bf{Theorem 2.3.}} {\it Let $a, b>0$ and $g\in L^2(K)$. If $\alpha_{0}, \gamma$ and $\sigma$ satisfy}

$$\sigma < \gamma\le \alpha_0 < \infty,\eqno(2.8) $$

\parindent=0mm \vspace{.1in}
{\it then the Gabor system ${\cal G}(g, a, b)$ given by (2.5)  constitutes a frame for $L^2(K)$ with bounds $\dfrac{C_{3}}{a}$ and $\dfrac{D_{3}}{a}$, where $C_{3}=\gamma-\sigma$ and $D_{3}= \alpha_{0}+\sigma.$ }

\parindent=0mm \vspace{.2in}
{\it Remark 2.} Since $\sigma \le \mu$, the frame bounds in Theorem 2.3 are better than ones in Theorem 2.2.

\parindent=0mm \vspace{.2in}
{\bf{3. Proof of the Main Results}}

\parindent=8mm \vspace{.1in}
In order to prove Theorems 2.2 and 2.3, we need the following lemma  whose proof can be found in [3].

\parindent=0mm \vspace{.1in}
{\bf {Lemma 3.1.}} {\it Suppose that $\left\{f_{k}\right\}_{k=1}^{\infty}$ is a family of elements in a Hilbert space $\mathbb H$ such that there exist constants $0<A\le B< \infty$ satisfying }

$$A\big\|f\big\|_{2}^{2}\le \sum_{k=1}^{\infty}\left|\big \langle f, f_{k}\big\rangle \right|^2\le B \big\|f\big\|_{2}^{2}, $$

\parindent=0mm \vspace{.1in}
{\it for all $f$ belonging to a dense subset ${\cal D} $ of $\mathbb H$. Then, the same inequalities are true for all $f\in \mathbb H$; that is, $\left\{f_{k}\right\}_{k=1}^{\infty}$  is a frame for $\mathbb H$. }

\parindent=8mm \vspace{.2in}
In view of Lemma 3.1, we will consider the following set of functions:

$$\Omega^0= \left\{f \in \Omega: \text{supp}\hat{f} \subset  K \backslash \left\{0\right\}~ \text{and}~\big\|\hat f\big\|_{\infty}<\infty\right\}.$$

\parindent=0mm \vspace{.1in}
Since $\Omega$ is dense in $L^2(K)$ and closed under the Fourier transforms, the set $\Omega^0$ is also dense in $L^2(\mathbb K)$. Therefore, it is enough to verify that the system ${\cal G}(g, a, b)$  given by (2.5) is a frame for $L^2(K)$ if the results of Theorems 2.2 and 2.3 hold for all $f \in \Omega^0$.

\parindent=8mm \vspace{.2in}
Assume that  $f\in L^2(K)$ and $h\in \Omega^{0}$, then by periodization, we have

$$\int_{K}h(\xi) f(\xi)\,\chi_{k}\big(a(\xi-\omega)\big)\, d\xi=\int_{G_{a^{-1}}} \sum_{k\in \mathbb N_{0}}h\big(\xi+a^{-1}u(k)\big) f\big(\xi+a^{-1}u(k)\big)\,\chi_{k}\big(a(\xi-\omega)\big)\, d\xi$$

\parindent=0mm \vspace{.1in}
Since $h$ lies in $\Omega^{0}$, so it is a bounded and compactly supported, therefore the number of $k$ in the above sum is finite. Thus, we can say the  series

$$\sum_{k\in \mathbb N_{0}}\int_{K}h(\xi) f(\xi)\,\chi_{k}\big(a(\xi-\omega)\big)\, d\xi\eqno(3.1)$$

\parindent=0mm \vspace{.1in}
is convergent to a periodic function $H(\xi)\in L^2(G_{a^{-1}})$, where
$$H(\xi)=\sum_{k\in \mathbb N_{0}}h\big(\xi+a^{-1}u(k)\big) f\big(\xi+a^{-1}u(k)\big),~ \text{and}~~G_{a}=\left\{ x\in K: |x|\le |a|\right\}.$$

\parindent=0mm \vspace{.1in}
{\bf{Proof of Theorem 2.2.}} For any $ f \in L^2(K)$, there exists a function sequence $\left\{ f_j\right\}_{j=1}^{\infty} \subset \Omega^{0}$, such that
$$\left\|\hat f_{j}- \hat f\right\|_{2} \to 0\quad \text{as}~j\to \infty,~\text{ and supp}~\hat f_{j} \subset {\mathfrak B}^{j}$$

\parindent=0mm \vspace{.1in}
since $\Omega^{0}$ is dense in $L^2(K)$. For fixed $m\in\mathbb N_{0}$, define the functional

$$ P_m(h)=\sum_{n \in \mathbb N_{0}}\left|\big \langle h, g_{m,n} \big \rangle \right|^2
= \sum_{n \in \mathbb N_{0}} \left|\big \langle \hat{h}, \hat{g}_{m,n}  \big \rangle \right|^2, \quad h \in L^2(K).\eqno(3.2)$$

\parindent=0mm \vspace{.1in}
Since the Fourier transform of

$$\hat g_{m,n}(\xi) =\overline{ \chi_{n}\Big(a\big(\xi-bu(m)\big)\Big)} \hat g\big(\xi-bu(m)\big),$$

\parindent=0mm \vspace{.1in}
therefore, by using equation (3.1), we are able to express (3.2) as
$$\begin{array}{rcl}
P_m(f_{j})&=&\displaystyle\sum_{n \in \mathbb N_{0}} \left|\left\langle \hat  f_{j}, \hat{g}_{m,n}\right\rangle \right|^2\\\\
&=&\displaystyle\sum_{n \in \mathbb N_{0}}\left\langle \hat  f_{j}, \hat{g}_{m,n}\right\rangle\overline{ \left\langle \hat  f_{j}, \hat{g}_{m,n}\right\rangle }\\\\
&=&\displaystyle\sum_{n \in \mathbb N_{0}}\int_{K}\hat f_{j}(\xi)\,\overline{\hat g\big(\xi-bu(m)\big)} \chi_{n}\Big(a\big(\xi-bu(m)\big)\Big)\,d\xi \\\ &&\qquad\qquad\qquad\times\displaystyle\int_{K}\overline{ \hat f_{j}(\omega)}\,\hat g\big(\omega-bu(m)\big) \overline{\chi_{n}\Big(a\big(\omega-bu(m)\big)\Big)}\,d\omega\\\\
&=&\displaystyle \dfrac{1}{a}\sum_{k\in \mathbb N_{0}}\int_{K} \hat f_{j}\big(\xi+a^{-1}u(k)\big)\overline{\hat g\big(\xi-bu(m)+a^{-1}u(k)\big)}\,\overline{ \hat f_{j}(\xi)}\,\hat g\big(\omega-bu(m)\big)\,d\xi.
\end{array}$$

Let
$$P(f) =\sum_{m\in \mathbb N_{0}}\sum_{n \in \mathbb N_{0}} \left|\big\langle f, g_{m,n}\big\rangle \right|^2= \sum_{m\in \mathbb N_{0}} P_{m}(f),\eqno(3.3)$$

then
\begin{align*}
P(f_{j}) &=\dfrac{1}{a}\sum_{k\in \mathbb N_{0}}\sum_{m\in \mathbb N_{0}}\int_{K}\hat f_{j}\big(\xi+a^{-1}u(k)\big)\overline{\hat g\big(\xi-bu(m)+a^{-1}u(k)\big)}\,\overline{ \hat f_{j}(\xi)}\,\hat g\big(\omega-bu(m)\big)\,d\xi \\\\
&= Q_{1}(f_j)+Q_{2}(f_j), \tag{3.4}
\end{align*}

\parindent=0mm \vspace{.0in}
where
\begin{align*}
 Q_{1}(f_j)&=\dfrac{1}{a}\sum_{m\in \mathbb N_{0}}\int_{K}\left|\hat f_{j}(\xi)\,\hat g\big(\omega-bu(m)\big) \right|^2 d\xi\tag{3.5}\\\\
  Q_{2}(f_j)&=\dfrac{1}{a}\sum_{k\in \mathbb N}\sum_{m\in \mathbb N_{0}}\int_{K}\hat f_{j}\big(\xi+a^{-1}u(k)\big)\overline{\hat g\big(\xi-bu(m)+a^{-1}u(k)\big)}\,\overline{ \hat f_{j}(\xi)}\,\hat g\big(\xi-bu(m)\big)\,d\xi.\tag{3.6}
\end{align*}

\parindent=0mm \vspace{.0in}
Since $\alpha_0 <\infty$, the series $Q_{1}(f_j)$ is convergent and

$$\dfrac{\gamma}{a}\left\|\hat f_{j}\right\|^{2}_{2}\le Q_{1}(f_j)\le \dfrac{\alpha_{0}}{a}\left\|\hat f_{j}\right\|^{2}_{2},$$

or equivalently
$$\dfrac{\gamma}{a}\,\big\|  f_{j}\big\|^2_{2}\le Q_{1}(f_j)\le \dfrac{\alpha_{0}}{a}\,\big\|  f_{j}\big\|^2_{2}.\eqno(3.7)$$

\parindent=0mm \vspace{.1in}
Next, we claim that $Q_{2}(f_j)$ is absolutely convergent. To prove this, we set

$$Q_{2}^{*}(f_j)=\dfrac{1}{a}\sum_{k\in \mathbb N}\sum_{m\in \mathbb N_{0}}\left|\int_{K}\hat f_{j}\big(\xi+a^{-1}u(k)\big)\overline{\hat g\big(\xi-bu(m)+a^{-1}u(k)\big)}\,\overline{ \hat f_{j}(\xi)}\,\hat g\big(\xi-bu(m)\big)\,d\xi\right|.$$

\parindent=0mm \vspace{.1in}
Note that
$$\Big| \hat{g}\big(\xi -bu(m)+a^{-1}u(k) \big) \hat{g}\big(\xi-bu(m) \big)\Big| \le \dfrac{1}{2}\left(\Big|\hat{g}\big( \xi-bu(m)+a^{-1}u(k)\big)\Big |^2 + \Big| \hat{g}\big(\xi-bu(m) \big)\Big|^2\right),$$

\parindent=0mm \vspace{.0in}
hence we have
$$Q_{2}^{*}(f_j)\le \dfrac{1}{a}\sum_{k\in \mathbb N}\sum_{m\in \mathbb N_{0}}\int_{K}\left|\hat f_{j}\big(\xi+bu(m)+a^{-1}u(k)\big) \overline{\hat f_{j}\big(\xi+a^{-1}u(k)\big)}\right|\big|\hat g(\xi)\big|^2\,d\xi.$$

\parindent=0mm \vspace{.1in}
Since each $f_{j}$ is bounded and compactly supported on ${\mathfrak B}^{j}$, and in fact they belongs to $\Omega^{0}$, hence there exist a constant $M>0$ such that

$$Q_{2}^{*}(f_j)\le M\left\|\hat f_{j}\right\|^2_{\infty}\big\|g\big\|^2 <\infty, $$

\parindent=0mm \vspace{.1in}
which proves our claim that $Q_{2}(f_j)$ is absolutely convergent.

\parindent=8mm \vspace{.1in}
Using Cauchy–Schwarz inequality, we obtain

\begin{align*}
\Big|Q_{2}(f_j)\Big|&=\left|\dfrac{1}{a}\sum_{k\in \mathbb N}\int_{K}\hat f_{j}\big(\xi+a^{-1}u(k)\big)\,\overline{ \hat f_{j}(\xi)}\,\Lambda_{k}(\xi)\,d\xi\right|\\\\
&\le \dfrac{1}{a}\sum_{k\in \mathbb N}\int_{K} \left\{\left| \hat f_{j}\big(\xi+a^{-1}u(k)\big)\right|\big|\Lambda_{k}(\xi)\big|^{1/2} \right\}\left\{\left| \hat f_{j}(\xi)\right|\big|\Lambda_{k}(\xi)\big|^{1/2} \right\}d\xi\\\\
&\le \dfrac{1}{a}\sum_{k\in \mathbb N}\left\{\int_{K} \left| \hat f_{j}\big(\xi+a^{-1}u(k)\big)\right|^2\big|\Lambda_{k}(\xi)\big|d\xi \right\}^{1/2}\left\{\int_{K}\left| \hat f_{j}(\xi)\right|^2\big|\Lambda_{k}(\xi)\big|d\xi \right\}^{1/2}.\tag{3.8}
\end{align*}

\parindent=0mm \vspace{.0in}
It is easy to verify that
$$\hat f_{j}\big(\xi+a^{-1}u(k)\big)= \hat f_{j}(\xi),\quad\text{and}\quad \Lambda_{k}\big(\xi-a^{-1}u(k)\big)=\Lambda_{k}(\xi), ~\forall~k\in\mathbb N.$$

\parindent=0mm \vspace{.1in}
Thus, we have
$$\Big|Q_{2}(f_j)\Big|\le \dfrac{1}{a}\int_{K}\left| \hat f_{j}(\xi)\right|^2 d\xi \sum_{k\in \mathbb N} \delta_{k}=\dfrac{\mu}{a}\,\big\|  f_{j}\big\|^2_{2},$$

or equivalently,

$$-\dfrac{\mu}{a}\,\big\|  f_{j}\big\|^2_{2}\le Q_{2}(f_j) \le \dfrac{\mu}{a}\,\big\|  f_{j}\big\|^2_{2}.\eqno(3.9)$$

\parindent=0mm \vspace{.1in}
It follows from (3.7) and (3.9) that
$$\dfrac{\gamma-\mu}{a}\,\big\|  f_{j}\big\|^2_{2}\le P(f_j) \le \dfrac{\alpha_{0}+\mu}{a}\,\big\|  f_{j}\big\|^2_{2}.$$

\parindent=0mm \vspace{.1in}
Letting $j\to \infty$ in above inequality, we obtain
$$\dfrac{\gamma-\mu}{a}\,\big\|  f\big\|^2_{2}\le P(f) \le \dfrac{\alpha_{0}+\mu}{a}\,\big\|  f\big\|^2_{2},$$

or
$$\dfrac{C_{2}}{a}\,\big\|  f\big\|^2_{2}\le \sum_{m\in \mathbb N_{0}}\sum_{n \in \mathbb N_{0}} \left|\big\langle f, g_{m,n}\big\rangle \right|^2 \le \dfrac{D_{2}}{a}\,\big\|  f\big\|^2_{2},$$

\parindent=0mm \vspace{.1in}
where $C_{2}=\gamma-\mu$ and $D_{2}={\alpha_{0}+\mu}.$ This completes the proof of Theorem 2.2.

\parindent=0mm \vspace{.2in}
{\bf{Proof of Theorem 2.3.}} Similar to the proof of Theorem 2.2, (3.5)--(3.8) hold. It follows from (3.8), the Cauchy-Schwarz inequality that

\begin{align*}
\Big|Q_{2}(f_j)\Big|&\le \dfrac{1}{a}\left\{\sum_{k\in \mathbb N}\int_{K} \left| \hat f_{j}\big(\xi+a^{-1}u(k)\big)\right|^2\big|\Lambda_{k}(\xi)\big|d\xi \right\}^{1/2}\left\{\sum_{k\in \mathbb N}\int_{K}\left| \hat f_{j}(\xi)\right|^2\big|\Lambda_{k}(\xi)\big|d\xi \right\}^{1/2}\\\\
&= \dfrac{1}{a}\left\{\sum_{k\in \mathbb N}\int_{K} \left| \hat f_{j}(\xi)\right|^2\big|\Lambda_{k}\big(\xi-a^{-1}u(k)\big)\big|d\xi \right\}^{1/2}\left\{\sum_{k\in \mathbb N}\int_{K}\left| \hat f_{j}(\xi)\right|^2\big|\Lambda_{k}(\xi)\big|d\xi \right\}^{1/2}\\\\
&= \dfrac{1}{a}\left\{\int_{K} \left| \hat f_{j}(\xi)\right|^2 \sum_{k\in \mathbb N}\big|\Lambda_{k}(\xi)\big|d\xi \right\}^{1/2}\left\{\int_{K}\left| \hat f_{j}(\xi)\right|^2\sum_{k\in \mathbb N} \big|\Lambda_{k}(\xi)\big|d\xi \right\}^{1/2}\\\\
&\le \dfrac{\sigma}{a}\,\big\|f_{j}\big|^2_{2},\tag{3.10}
\end{align*}

\parindent=0mm \vspace{.1in}
which implies that

$$-\dfrac{\sigma}{a}\,\big\|  f_{j}\big\|^2_{2}\le Q_{2}(f_j) \le \dfrac{\sigma}{a}\,\big\|  f_{j}\big\|^2_{2}.\eqno(3.11)$$

\parindent=0mm \vspace{.1in}
Combining (3.7) with (3.11), we obtain
$$\dfrac{\gamma-\sigma}{a}\,\big\|  f_{j}\big\|^2_{2}\le P(f_j) \le \dfrac{\alpha_{0}+\sigma}{a}\,\big\|  f_{j}\big\|^2_{2}.$$

\parindent=0mm \vspace{.1in}
By taking $j\to \infty$ in above relation, we get
$$\dfrac{\gamma-\sigma}{a}\,\big\|  f\big\|^2_{2}\le P(f) \le \dfrac{\alpha_{0}+\sigma}{a}\,\big\|  f\big\|^2_{2},$$

or
$$\dfrac{C_{3}}{a}\,\big\|  f\big\|^2_{2}\le \sum_{m\in \mathbb N_{0}}\sum_{n \in \mathbb N_{0}} \left|\big\langle f, g_{m,n}\big\rangle \right|^2 \le \dfrac{D_{3}}{a}\,\big\|  f\big\|^2_{2},$$

\parindent=0mm \vspace{.1in}
where $C_{3}=\gamma-\sigma$ and $D_{3}={\alpha_{0}+\sigma}.$ This completes the proof of Theorem 2.3.

\parindent=0mm \vspace{.2in}

{\bf{References}}


\begin{enumerate}

{\small {

\bibitem[1]{1} R. J. Duffin and A. C. Shaeffer, A class of nonharmonic Fourier series, {\it Transactions of the American
Mathematical Society}, vol. 72, pp. 341-366, 1952.

\bibitem[2]{2}  I. Daubechies, A. Grossmann and Y. Meyer, Painless non-orthogonal expansions,  {\it Journal of Mathematical Physics},  vol.  27, no. 5, pp. 1271-1283, 1986.

\bibitem[3]{3} O. Christensen, {\it An Introduction to Frames and Riesz Bases}, Birkh\"{a}user, Boston, 2015.

\bibitem[4]{4} L. Debnath and F. A. Shah, {\it Wavelet  Transforms and Their Applications},  Birkh\"{a}user, New York, 2015.


\bibitem[5]{5}  A. Ron and Z. Shen,  Weyl-Heisenberg frames and Riesz bases in $L^2(\mathbb R^d )$, {\it  Duke Mathematics Journal}, vol. 89, pp. 237-
282, 1997.

\bibitem[6]{6} K. Gr\"{o}chenig,  A. J. Janssen, N. Kaiblinger and GE. Pfander,  Note on $B$-splines, wavelet scaling functions, and Gabor frames, {\it IEEE Transactions and Information Theory}, vol. 49, no. 12, pp. 3318-3320, 2003.

\bibitem[7]{7} P. G. Casazza and O. Christensen, Weyl-Heisenberg frames for subspaces of $L^2(\mathbb R)$,  {\it Proceedings of American Mathematical Society}, vol. 129, pp. 145-154, 2001.

\bibitem[8]{8} K. Wang, Necessary and sufficient conditions for expansions of Gabor type, {\it Analysis in Theory and Applications}, vol. 22, pp. 155-171,  2006.

\bibitem[9]{9} X. L. Shi and F. Chen, Necessary conditions for Gabor frames, {\it Science in China : Series A.} vol. 50, no. 2, pp. 276-284, 2007.

\bibitem[10]{10} D. Li, G. Wu  and X. Zhang, Two sufficient conditions in frequency domain for Gabor frames, {\it Applied Mathematics Letters}, vol. 24, pp. 506-511, 2011.

\bibitem[11]{11} K. Gr\"{o}chenig, {\it Foundation of Time-Frequency Analysis,} Birkh\"{a}user, Boston, 2001.

\bibitem[12]{12} H. G. Feichtinger and T. Strohmer, {\it Advances in Gabor Analysis}, Birkh\"{a}user, Boston, 2003.

\bibitem[13]{13} M. H. Taibleson, {\it Fourier Analysis on Local Fields}, Princeton University Press, Princeton, NJ, 1975.

\bibitem[14]{14} D. Li and H. K. Jiang,  Basic results Gabor frame on local fields, {\it Chinese Annals of Mathematics: Series B}, vol. 28, no. 2, 165-176, 2007.

\bibitem[15]{15} F. A. Shah, A characterization of tight Gabor frame on local fields of positive characteristic, {\it Journal of the Nigerian Mathematical Society}, (accepted), 2015.

\bibitem[16]{16} F. A. Shah, Gabor frames on a half-line, {\it Journal of Contemporary Mathematical Analysis}, vol. 47, no. 5, pp. 251-260, 2012.

\bibitem[17]{17} F. A. Shah, Frame multiresolution analysis on local fields of positive characteristic, {\it Journal of Operators}. Article ID 216060, 8 pages, 2015.

\bibitem[18]{18}  F. A. Shah and Abdullah, Wave packet frames on local fields of positive characteristic,  {\it Applied Mathematics and Computation,} vol.  249, pp. 133-141, 2014.

\bibitem[19]{19} F. A. Shah and Abdullah, A characterization of tight wavelet frames on local fields of positive characteristic, {\it Journal of Contemporary Mathematical Analysis}, vol. 49, pp. 251-259, 2014.

\bibitem[20]{20} F. A. Shah and L. Debnath, Tight wavelet frames on local fields, {\it Analysis}, vol. 33, pp. 293-307, 2013.
}}


\end{enumerate}
\end{document}